\documentclass[aos,nameyear,seceqn,dvips]{arximspdf}

\doi{10.1214/07-AOS543}
\volume{36}
\issue{5}
\pubyear{2008}
\firstpage{2453}
\lastpage{2470}

\makeatletter
\newproclaim{assumption}{Assumption}[section]
\newtheorem{theorem}{Theorem}[section]
\newproclaim{remark}{Remark}[section]
\newtheorem{corollary}{Corollary}[section]
\newtheorem{lemma}{Lemma}[section]
\makeatother

\begin{document}
\begin{frontmatter}

\title{Residual empirical processes for long and short memory
time series\protect\thanksref{T1}}
\runtitle{Residual empirical processes}
\thankstext{T1}{Supported in part by CERG from the Hong Kong Research
Grants Council Grant Numbers:
CUHK400305, CUHK400306, HKUST6022/05P, HKUST6428/06H and HKUST6016/07.}

\begin{aug}
\author[A]{\fnms{Ngai Hang} \snm{Chan}} \and
\author[B]{\fnms{Shiqing} \snm{Ling}\corref{}\ead[label=e2]{maling@ust.hk}}
\runauthor{N. H. Chan and S. Ling}
\affiliation{Chinese University of Hong Kong and\break Hong Kong
University of Science and Technology}
\address[A]{Department of Statistics\\
Chinese University of Hong Kong,\\
Shatin, NT\\
Hong Kong} 
\address[B]{Department of Mathematics\\
Hong Kong University of Science \\
\quad and Technology\\
Hong Kong\\
\printead{e2}}
\end{aug}

\received{\smonth{6} \syear{2006}}
\revised{\smonth{8} \syear{2007}}

%
\begin{abstract}
This paper studies the residual empirical process of long- and
short-memory time series regression models and establishes its
uniform expansion under a general framework. The results are
applied to the stochastic regression models and unstable
autoregressive models. For the long-memory noise, it is shown that
the limit distribution of the Kolmogorov--Smirnov test statistic
studied in Ho and Hsing [\textit{Ann. Statist.} \textbf{24} (1996) 992--1024]
does not hold when the stochastic
regression model includes an unknown intercept or when the
characteristic polynomial of the unstable autoregressive model has a
unit root. To this end, two new statistics are proposed to test for
the distribution of the long-memory noises of stochastic regression
models and unstable autoregressive models.
\end{abstract}

%
\begin{keyword}[class=AMS]
\kwd[Primary ]{62G30}
\kwd[; secondary ]{62M10}.
\end{keyword}
\begin{keyword}
\kwd{Empirical process}
\kwd{long-memory time series}
\kwd{residuals}
\kwd{unit root}
\kwd{weak convergence}.
\end{keyword}

\end{frontmatter}

\section{Introduction}\label{s1}

Let the time
series $\{y_{t}\}$ be generated by the model
%
\begin{eqnarray}\label{e1.1}
y_{t}=\beta'X_{t}+ \varepsilon_{t}
\quad \mbox{and}\quad
\varepsilon_{t}=\sum_{i=0}^{\infty}a_{i}e_{t-i} ,
\end{eqnarray}
where $X_{t}$'s are a sequence of $p$-dimensional time series which
are measurable with respect to
$\mathcal{F}_{t-1}=\sigma\{\varepsilon_{t-1}, \varepsilon_{t-2},\ldots\}$ or
independent of $\{\varepsilon_{t}\}$. The coefficients $a_i$ satisfy
$\sum_{i=1}^{\infty} a_{i}^{2}<\infty$; $a_{0}=1$ and
$a_{k}=k^{H-3/2}L_{0}(k)$ for some slowly varying function $L_{0}$
[see Feller (\citeyear{fel1971})] with $H<1$; and $\{e_{t}\}$ is a sequence of i.i.d.
mean zero random variables with $\sigma_{e}^{2}=Ee_{t}^{2} <
\infty$. The process $\{\varepsilon_{t}\}$ exhibits a long-memory
(short-memory) phenomenon when $H\in(1/2,1)$ ($H<1/2$), which has
been considerably studied in the literature; see, for example,
Robinson (\citeyear{rob1995a}, \citeyear{rob1995b})
and the references therein. When
model (\ref{e1.1}) is used to construct forecasting intervals or
value-at-risk (VaR), knowledge on the distribution function $F(x)$
of $\varepsilon_{t}$ is of crucial importance. This motivates the
study on testing of $F(x)$ and
on related empirical processes of $\{\varepsilon_{t}\}$.

When $H\in(1/2,1)$, Ho and Hsing (\citeyear{hoHsi1996})
established a strong expansion for the
empirical process of $\{\varepsilon_{t}\}$ in (\ref{e1.1}). Specifically,
let
%
\begin{equation}\label{e1.2}
K_{n}(x)=\frac{1}{\sigma_{n}}\sum_{t=1}^{n}[I(\varepsilon_{t}
\le x)-F(x)] ,
\end{equation}
where $I(\cdot)$ is the indicator function and
$\sigma_{n}^{2}=\operatorname{var}(\sum_{t=1}^{n}\varepsilon_{t})$.
They proved that
%
\begin{eqnarray}
&\displaystyle \sup_{x}\Biggl|K_{n}(x)+\frac{1}{\sigma_{n}}F'(x)
\sum_{t=1}^{n}\varepsilon_{t}\Biggr|=o(1) \qquad \mbox{a.s.}, &
\label{e1.3}\\
&\displaystyle
\sigma_{n}^{2}\sim\kappa(H)n^{2H}L_{0}^{2}(n)
\quad \mbox{and} \quad \sigma_{n}^{-1}\sum_{t=1}^{n}
\varepsilon_{t}\stackrel{\mathcal{L}}{\to} N(0,1) ;
\label{e1.4}
\end{eqnarray}
see also Taqqu (\citeyear{taq1975}) and Hosking (\citeyear{hos1996}). Herein,
$\sup_{x}=\sup_{x\in R}$,
$\kappa(H)=\int_{0}^{\infty}(x+x^{2})^{H-3/2} \,dx$, $a_{n}\sim b_{n}$
means that $a_{n}/b_{n}\to1$ as $n\to\infty$ and $\stackrel{\mathcal{L}}{\to}$
denotes convergence in distribution as $n\to\infty$. By (\ref{e1.3}),
%
\begin{equation}\label{e1.5}
\biggl[
\sup_{x}F'(x)\biggr]^{-1}
\sup_{x}|K_{n}(x)|\stackrel{\mathcal{L}}{\to} | N(0,1)| ,
\end{equation}
if $\sup_{x}|F'(x)|<\infty$. This is the Kolmogorov--Smirnov test
statistic of Ho and Hsing (\citeyear{hoHsi1996}) for testing the distribution
$F(x)$. Contrary to the standard weak convergence of the empirical
process in the short-memory case, the result (\ref{e1.5}) is somewhat
striking as $\sup_{x}|K_{n}(x)|$ does not converge to the maximum of
a Brownian bridge as in the traditional case. Weak convergence of
$\{K_{n}(x)\}$ was established in Dehling and
Taqqu (\citeyear{dehTaq1989}) when $\{\varepsilon_{t}\}$ is a long-range dependent
Gaussian process. Koul and Surgailis (\citeyear{kuoSur1997}) obtained some related
results when $H\in(1/2,1)$. Wu (\citeyear{wu2003}) showed that (\ref{e1.3}) holds in
probability under a weaker condition and a general setup and
characterized the limit behavior of $K_{n}(x)$ when $H \le1/2$;
see also Ho and Hsing (\citeyear{hoHsi1997}).

Note that since $\{\varepsilon_{t}\}$ is unobservable in model
(\ref{e1.1}), the Kolmogorov--Smirnov test has to be evaluated based on the
residual process of $\{\varepsilon_{t}\}$. In this situation, a key
issue of interest is to determine the validity of (\ref{e1.5}) for the
Kolmogorov--Smirnov statistic when $\{\varepsilon_{t}\}$ is replaced
by its corresponding residual process. Furthermore, when (\ref{e1.5}) becomes
invalid, how can one test for the
distribution of $\{\varepsilon_{t}\}$? These two issues have been
studied extensively when $\{\varepsilon_{t}\}$ is i.i.d.; see Bai
(\citeyear{bai1994}, \citeyear{bai1996}, \citeyear{bai2003}),
Ling (\citeyear{lin1998}),
Lee and Wei (\citeyear{leeWei1999}),
Koul (\citeyear{kou2002}),
Lee and Taniguchi (\citeyear{leeTan2005}) and
Koul and Ling (\citeyear{kouLin2006}) for further discussions.
But for model (\ref{e1.1}) and for the Kolmogorov--Smirnov statistic studied in
Ho and Hsing (\citeyear{hoHsi1996}),
these two important issues still remain unresolved. When
$\beta'X_{t}$ is a constant and $\varepsilon_{t}$ is an
ARFIMA$(p,d,q)$ model, the distribution of $\{\varepsilon_t\}$ can be
determined by $\{e_t\}$ once the parameters of the ARFIMA model are
estimated. In this case, it would be sufficient to test for the
distribution of $\{e_{t}\}$, for which standard
procedures for residuals from a model with i.i.d. noises, such as those
given in
Bai (\citeyear{bai1994}) and
Lee and Wei (\citeyear{leeWei1999}), can be adopted.
To study the general residual process of $\{\varepsilon_{t}\}$,
however, substantially different arguments need to
be employed which rely heavily on the results of
Ho and Hsing (\citeyear{hoHsi1996}, \citeyear{hoHsi1997})
and Wu (\citeyear{wu2003}).

This paper first establishes a uniform expansion of the residual
empirical process of $\{\varepsilon_t\}$
under a general framework. The result is used to study
the stochastic regression model of Robinson and Hidalgo (\citeyear{robHid1997}) and
the unstable AR model of Chan and Terrin (\citeyear{chaTer1995}),
Truong-Van and Larramendy (\citeyear{troLar1996})
and Wu (\citeyear{wu2006}). It is shown
that the test statistic (\ref{e1.5}) of Ho and Hsing (\citeyear{hoHsi1996}) is no longer
valid when the stochastic regression model includes an unknown
intercept or when the characteristic
polynomial of the unstable AR model has a unit root. Our results not only
encompass the long-memory $\{\varepsilon_{t}\}$, but also the
short-memory $\{\varepsilon_{t}\}$. Furthermore,
two new statistics are constructed to test the distribution of
the long-memory noises in the stochastic regression model and the
unstable AR model.

This paper is
organized as follows. A general result is given in Section \ref{s2}.
The residual processes of stochastic regression and unstable time
series are
presented in Sections \ref{s3} and \ref{s4}, respectively.

\section{A general result}\label{s2}

Let $\hat{\beta}_{n}$ be an estimator of $\beta$ in (\ref{e1.1}). Let
$\hat{\varepsilon}_{t}= y_{t}-\hat{\beta}_{n}'X_{t}$ be the residual
of model (\ref{e1.1}). Further, define the empirical process based on
residuals $\{\hat{\varepsilon}_{t}\}$ by
\begin{eqnarray*}
\hat{K}_{n}^{\delta}(x)=\frac{1}{\sigma_{n}}\sum_{t=1}^{n}
[I(\hat{\varepsilon}_{t}\le x)-F(x)] .
\end{eqnarray*}
For $H\in(1/2, 1)$, $\sigma_{n}$ is given in (\ref{e1.4}). For
$\sum_{j=0}^{\infty}|a_{j}|<\infty$, which implies \mbox{$H\le1/2$}, Ho
and Hsing (\citeyear{hoHsi1997}) show that
$\sigma^{2}\equiv\lim_{n\to\infty}\sigma_{n}^{2}/n$ exists and is
finite; see also Wu (\citeyear{wu2003}).
Let $G_{0}$ be the common
distribution of $\{e_{t}\}$. Write $\varepsilon_{t}=e_{t}+\xi_{t-1}$
and let $A_{t}(x)=G'_{0}(x-\xi_{t-1})-E[G'_{0}(x-\xi_{t-1})]$, where
$\xi_{t-1}=\sum_{i=1}^{\infty}a_{i}e_{t-i}$. Denote
$\|\cdot\|=\operatorname{tr}(M'M)$ for some matrix or vector $M$.
We need the following two assumptions.
\begin{assumption}\label{as2.1}
(a) $H < 1/2$ and $\sigma> 0$, or
$H=1/2$, $\sigma> 0$ and $\sum_{j=0}^{\infty}|a_{j}|<\infty$, or
$1/2 < H < 1$, and (b) $G_{0}$ is three times
differentiable with bounded, continuous and integrable derivatives
such that $\int x^{4} \,dG_{0}(x)<\infty$.
\end{assumption}
\begin{assumption}\label{as2.2}
Let $\delta_{n}$ be a $p\times p$ constant
matrix depending on $n$ such that the following statements hold:
\begin{longlist}[d]
\item[(a)] $\delta_{n}^{-1}(\hat{\beta}_{n}-\beta)=O_{p}(1)$,
\item[(b)] $\sigma_{n}^{-1}\sum_{t=1}^{n} E\|\delta_{n}'X_{t}\|=O(1)$,
\item[(c)] $\sigma_{n}^{-1}\sum_{t=1}^{n}E\|\delta_{n}'X_{t}\|^{2}=o(1)$,
\item[(d)] $\sigma_{n}^{-1}\sup_{x}\|\sum_{t=1}^{n}A_{t}(x)
\delta_{n}'X_{t}\|=o_{p}(1)$.
\end{longlist}
\end{assumption}

Assumption \ref{as2.1}(b) can be replaced by a general condition
in Wu (\citeyear{wu2003}).
$\delta_{n}$~is the rate of convergence of $\hat{\beta}_{n}$.
Assumptions \ref{as2.2}(b) and (c) automatically hold if
$\delta_{n}^{-1}=\sqrt{n}I_{p}$ and $X_{t}$ is strictly stationary
with $E\|X_{t}\|^{2}<\infty$, where $I_{p}$ is the \mbox{$p\times p$}
identity matrix. As will be seen in Sections \ref{s3} and \ref{s4},
$\delta_{n}^{-1}$ may not always be equal to $\sqrt{n}I_{p}$.
Assumptions \ref{as2.2}(b)--(d) are sufficient for the remainder term in the
following expansion to be negligible, although they may not be the
weakest ones. We state a general result as follows.
\begin{theorem}\label{thm2.1}
Assume that Assumption~\ref{as2.1} and
Assumption~\ref{as2.2} hold. Then
\begin{eqnarray*}
\sup_{ x }|\hat{K}_{n}( x)-K_{n}( x)- R_{n}F'(x)|=o_{p}(1) ,
\end{eqnarray*}
where $R_{n}=\sigma_{n}^{-1}(\hat{\beta}_{n}-\beta)'\sum_{t=1}^{n}X_{t}=O_{p}(1)$.
\end{theorem}
\begin{remark}\label{rem2.1}
According to this theorem, if $R_n=o_p(1)$,
then $\sup_{ x}|\hat{K}_{n}( x)-K_{n}( x)|=o_{p}(1)$ and, hence,
$\sup_x |\hat{K}_n(x)|$ and $\sup_x |K_n(x)|$ have the same
limit distribution. If $R_n \neq o_p(1)$,
then the limit distribution of $\sup_x |\hat{K}_n(x)|$ may be different
from that of $\sup_x |K_n(x)|$, as seen in Theorems \ref{thm3.1}
and \ref{thm4.1}. When $H\in(1/2,1)$, $K_{n}( x)$ can be replaced
by $-F'(x)\sum_{t=1}^{n} \varepsilon_{t}/\sigma_{n}$. When
$H<1/2$ with $EX_{t}=0$ or when $H \in(1/2,1)$, 
$\delta_{n}^{-1}=\sqrt{n}I_{p}$ and 
$\{X_{t}\}$ is strictly stationary, then $R_{n}=o_{p}(1)$.
\end{remark}
\begin{remark}\label{rem2.2}
We require $\{a_{k}\}$ to have the form
$k^{H-3/2}L_{0}(k)$ because we have to use the tightness condition
of empirical processes of $\{\varepsilon_{t}\}$ of Ho and Hsing (\citeyear{hoHsi1996})
and Wu (\citeyear{wu2003}) for $H\in(1/2, 1)$; and Theorem 3 and Corollary 2 of Wu
(\citeyear{wu2003}) for $H\le1/2$. Without this
condition, Theorem~\ref{thm2.1} is still valid
if $\sum_{i=0}^{\infty}|a_{i}|<\infty$ as long as the empirical
process of $\{\varepsilon_{t}\}$ is tight on $R$.
\end{remark}
\begin{pf*}{Proof of Theorem \ref{thm2.1}}
Let $\hat{u}_{n}=\delta^{-1}_{n}(\hat{\beta}_{n}-\beta)$. Then
$\hat{\varepsilon}_{t}=\varepsilon_{t}-\hat{u}_{n}'\delta_{n}'X_{t}$ and
\begin{eqnarray*}
&& \hat{K}_{n}(x)-K_{n}(x)-\frac{1}{\sigma_{n}}\sum_{t=1}^{n}
F'(x)\hat{u}_{n}'\delta_{n}'X_{t}
\\
&&\qquad=\frac{1}{\sigma_{n}}\sum_{t=1}^{n}[I(\varepsilon_{t}\le x
+ \hat{u}'_{n}\delta_{n}'X_{t})-I(\varepsilon_{t}\le x)
- F'(x)\hat{u}_{n}'\delta_{n}'X_{t}] .
\end{eqnarray*}
To study the process $\hat{K}_{n}(x)$, consider the process
\begin{eqnarray*}
A_{n}(x, u)=\frac{1}{\sigma_{n}}\sum_{t=1}^{n}
[I(\varepsilon_{t}\le x+u'\delta_{n}'X_{t})-I(\varepsilon_{t}\le x)-u'
F'(x)\delta_{n}'X_{t}]
\end{eqnarray*}
for all $u\in R^{p}$ and $x\in R$. By Assumption~\ref{as2.2}(a), if we can
show that
%
\begin{eqnarray}\label{e2.1}
\sup_{u\in[-\Delta, \Delta]^{p}}\sup_{ x}|A_{n}(x, u)|=o_{p}(1)
 \qquad   \mbox{for every }     \Delta\in(0, \infty),
\end{eqnarray}
then Theorem~\ref{thm2.1} is proved.
Denote
\begin{eqnarray*}
Z_{n}(x, u)=\frac{1}{\sigma_{n}}\sum_{t=1}^{n}[I(
\varepsilon_{t}\le x+u'\delta_{n}'X_{t})-F(x+u'\delta_{n}'X_{t})-
I(\varepsilon_{t}\le x)+F(x)].
\end{eqnarray*}
By the triangular inequality, $|A_{n}(x, u)|\le|Z_{n}(x,
u)|+|H_{n}(x, u)|$, where
\begin{eqnarray*}
H_{n}(x, u)=\frac{1}{\sigma_{n}}\sum_{t=1}^{n}[F(x+u'\delta_{n}'X_{t})-F(x)
-u'\delta_{n}'X_{t}F'(x)] .
\end{eqnarray*}
Since $\sup_{x}|G_{0}''(x)|<\infty$, we have
$\sup_{x}|F''(x)|<\infty$. Using this fact, Assumption~\ref{as2.2}(c) and
the Taylor expansion, $\sup_{u\in[-\Delta, \Delta]^{p}}\sup_{
x}|H_{n}(x, u)|=o_{p}(1)$. To prove (\ref{e2.1}), it is sufficient to show
that the following equation holds:
%
\begin{equation}\label{e2.2}
\sup_{u\in[-\Delta, \Delta]^{p}}\sup_{ x}|Z_{n}(x, u)|=o_{p}(1),
\end{equation}
for every $\Delta>0$.
For each $u\in R^{p}$ and $\lambda\in R$, let
%
\begin{eqnarray}\label{e2.3}
\tilde{Z}_{n}(x,u,\lambda) &=&
\frac{1}{\sigma_{n}}\sum_{t=1}^{n}\bigl[I\bigl(
\varepsilon_{t}\le x+g_{t}(u, \lambda)\bigr)
\nonumber\\[-8pt]
\\[-8pt]
&&\hspace*{33pt}{}
- F\bigl(x+g_{t}(u, \lambda)\bigr)
- I(\varepsilon_{t}\le x)+F(x)\bigr],
\nonumber
\end{eqnarray}
where $g_{t}(u, \lambda)= u'\delta_{n}'X_{t}+\lambda\|\delta
_{n}'X_{t}\|$. For every $\delta>0$, partition the rectangle $[-\Delta,
\Delta]^{p}$ into $m$ balls $\{C_{1}, \ldots, C_{m}\}$ each with
radius $\delta$. Take one point in each $C_{r}$ and denote it by
$u_{r}$. For any $u \in C_{r}$, we have
%
\begin{equation}\label{e2.4}
|g_{t}(u, \lambda)-g_{t}(u_{r}, \lambda)|\le
\|u-u_{r}\|\|\delta_{n}'X_{t}\| \le\delta\|\delta_{n}'X_{t}\| .
\end{equation}
Thus, $g_{t}(u_{r}, \lambda-\delta)\le g_{t}(u, \lambda)\le
g_{t}(u_{r}, \lambda+\delta)$. Note that $Z_{n}(x,
u)=\tilde{Z}_{n}(x, u, 0)$. By the monotonicity of the indicator
function, we obtain that
%
\begin{equation}\label{e2.5}
\qquad
Z_{n}(x, u)\le\tilde{Z}_{n}(x, u_{r},
\delta)+\frac{1}{\sigma_{n}}\sum_{t=1}^{n}\bigl[F\bigl(x+g_{t}(u_{r},
\delta)\bigr)-F\bigl(x+g_{t}(u, 0)\bigr)\bigr]
\end{equation}
and a reverse inequality holds when $\delta$ is replaced by
$-\delta$. Since $\sup_{x}|G'_{0}(x)|<\infty$, we
have $\sup_{x}|F'(x)|<\infty$. By the mean value theorem,
when $u\in C_{r}$,
%
\begin{eqnarray}\label{e2.6}
&& \Biggl |\frac{1}{\sigma_{n}}\sum_{t=1}^{n}
\bigl[F\bigl(x+g_{t}(u_{r},\pm \delta)\bigr)
- F\bigl(x+g_{t}(u, 0)\bigr)\bigr] \Biggr|
\nonumber\\
&&\qquad \le\frac{\sup_{x}|F'(x)|}{\sigma_{n}}
\sum_{t=1}^{n}|g_{t}(u_{r}, \pm \delta)-g_{t}(u, 0)|
\\
&&\qquad \le\frac{O(1)\delta}{\sigma_{n}}
\sum_{t=1}^{n}\|\delta_{n}'X_{t}\|= O_{p}(\delta) ,
\nonumber
\end{eqnarray}
where the last equality follows from Assumption~\ref{as2.2}(b) and the $O_{p}(1)$
holds uniformly for all $x\in\tilde{R}$, all $u\in C_{r}$ and all
$r=1,\ldots,m$.

Given any $\varepsilon>0$ and $\eta>0$, by (\ref{e2.6}), there
exists a $\delta_{1\varepsilon}>0$ such that
\begin{eqnarray*}
P\Biggl\{\frac{1}{\sigma_{n}}\max_{r}
\max_{u\in C_{r}}\sup_{x}\Biggl|\sum_{t=1}^{n}
\bigl[F\bigl(x+g_{t}(u_{r}, \pm \delta)\bigr)
- F\bigl(x+g_{t}(u, 0)\bigr)\bigr]\Biggr|
\ge\frac{\varepsilon}{3}\Biggr\}
\le\frac{\eta}{6} ,
\end{eqnarray*}
when $\delta\le\delta_{1\varepsilon}$ and $n\rightarrow\infty$.
By Lemma \ref{lemA.3}, there exists a $\delta_{2\varepsilon}>0$ such that
\begin{eqnarray*}
P\biggl\{\max_{r}\sup_{ x}|\tilde{Z}_{n}(x, u_{r}, \pm\delta)|\ge
\frac{\varepsilon}{3}\biggr\}
&\le&  P\biggl\{\max_{r}J_{3n}(u_{r}, \pm\delta)\ge
\frac{\varepsilon}{6}\biggr\}
+ P\biggl\{ \delta J_{4n}\ge \frac{\varepsilon}{6}\biggr\}
\\
&\le&  m \max_{r}P\biggl\{J_{3n}(u_{r}, \pm\delta)\ge
\frac{\varepsilon}{6}\biggr\}
+ \frac{\eta}{6}\le\frac{\eta}{3},
\end{eqnarray*}
when $\delta\le\delta_{2\varepsilon}$ and $n\rightarrow\infty$
because $m$ is an integer depending on $\delta$ but not depending on
$n$. By the preceding two inequalities, when $\delta\le
\min\{\delta_{1\varepsilon}, \delta_{1\varepsilon}\}$,
\begin{eqnarray*}
&& P\biggl\{\sup_{u\in[-\Delta, \Delta]^{p}}\sup_{ x}|Z_{n}(x,
u)| \ge\varepsilon\biggr\}
\\
&&\qquad \le P\biggl\{\max_{r}\sup_{ x}|\tilde{Z}_{n}(x, u_{r},
\delta)|\ge\frac{\varepsilon}{3}\biggr\}
+ P\biggl\{\max_{r}\sup_{x}|\tilde{Z}_{n}(x, u_{r}, -\delta)|
\ge\frac{\varepsilon}{3}\biggr\}
\\
&&\qquad \quad{} + P\Biggl\{\frac{1}{\sigma_{n}}
\max_{r}\max_{u\in C_{r}}\sup_{x}
\Biggl|\sum_{t=1}^{n} \bigl[F\bigl(x+g_{t}(u_{r}, \pm\delta)\bigr)
- F\bigl(x+g_{t}(u, 0)\bigr)\bigr]\Biggr|
\ge\frac{\varepsilon}{3}\Biggr\}
\\
&&\qquad\le\eta , \qquad \mbox{when } n \rightarrow\infty, \mbox{
proving (\ref{e2.2}).}
\end{eqnarray*}\upqed
\end{pf*}

\section{Residual empirical process of stochastic regression
models}\label{s3}

In this section we apply the results
in Section \ref{s2} to the stochastic regression model of Robinson and
Hidalgo (\citeyear{robHid1997}):
%
\begin{equation}\label{e3.1}
y_{t}=\alpha_{0}+\alpha' x_{t}+\varepsilon_{t},
\end{equation}
where $\varepsilon_{t}$ is defined in model (\ref{e1.1}), $x_{t}$ is a
$q$-dimension vector time series independent of $\{\varepsilon_{t}\}$,
and $\beta=(\alpha_{0}, \alpha')'$ is a $p=q+1$
dimensional unknown parameter vector. The least squares
estimator (LSE) or generalized LSE of $\alpha$ is not
asymptotically normal when both $x_{t}$ and $\varepsilon_{t}$
exhibit long-range dependence; see Robinson (\citeyear{rob1994}).
Robinson and Hidalgo (\citeyear{robHid1997})
proposed a class of weighted LSE which is $\sqrt{n}$-consistent
and asymptotically normal.

Let $f(\lambda)$ be the spectral density of $\varepsilon_{t}$ and $\phi
(\lambda)$ be a real-valued, even and integrable periodic
function with period $2\pi$ such that
$\psi(\lambda)=\phi^{2}(\lambda)f(\lambda)$ is continuous. Denote
$\phi_{j}=(2\pi)^{-2}\int_{-\pi}^{\pi}\phi(\lambda)\cos j\lambda
\, d\lambda$. Robinson--Hidalgo's weight\-ed LSE of
$\alpha$ is defined as
\[
\hat{\alpha}_{n}= \Biggl[\sum_{t=1}^{n}\sum_{s=1}^{n}(x_{t}-\bar{x})(x_{s}-\bar{x})'
\phi_{t-s} \Biggr]^{-1}
\Biggl [\sum_{t=1}^{n}\sum_{s=1}^{n}(x_{t}-\bar{x})(y_{s}-\bar{y})\phi_{t-s} \Biggr] ,
\]
where $\bar{x}=\sum_{t=1}^{n}x_{t}/n$ and
$\bar{y}=\sum_{t=1}^{n}y_{t}/n$. Let $\gamma_{j}=E(\varepsilon
_{t}\varepsilon_{t+j})$ and $\kappa_{abcd}(s,u,v,\break w)$ be
the fourth cumulant of $x_{as}$, $x_{bu}$, $x_{cv}$ and $x_{dw}$,
where $x_{as}$ is the $a$th element of $x_{s}$. Recall the assumptions
of Robinson and Hidalgo (\citeyear{robHid1997}) as follows.
\begin{assumption}\label{as3.1}
(a) $\sum_{j=0}^{\infty}\tilde{\phi}_{j}<\infty$ and
$(\sum_{j=0}^{n}|\gamma_{j}|+n\tilde{\gamma}_{n})[(\sum_{j=0}^{n}\tilde
{\phi}_{j}^{1/2})^{2}+n\Phi_{n}]=O(n)$ as
$n\to\infty$, where $\tilde{\gamma}_{a}=\max_{j\ge a}|\gamma_{j}|$,
$\tilde{\phi}_{a}=\max_{j\ge a}|\phi_{j}|$ and
$\Phi_{a}=\sum_{|j|>a}|\phi_{j}|$.

(b) $\{x_{t}\}$ is fourth-order stationary,
$\Gamma_{u}=E[(x_{1}-Ex_{1})(x_{1+|u|}-Ex_{1})']\to0$ and
$\max_{|v|,|w|<\infty}|\kappa_{abcd}(0,u,v,w)|\to0$ as
$|u|\to\infty$, $1\le a,b,c,d\le q$.

(c) $\Sigma_{\psi}$ is finite and $\Sigma_{\phi}$ and $\Sigma_{\psi}$
are nonsingular, where
$\Sigma_{\chi}=\int_{-\pi}^{\pi}\chi(\lambda)  \,dH(\lambda)/\break(2\pi)$
and $H(\lambda)$ is the Hermitian matrix such that $\Gamma_{j}=
\int_{-\pi}^{\pi}e^{ij\lambda}  \, dH(\lambda)$.
\end{assumption}

Discussions on this assumption, the choice of $\phi$ and its
computational procedures can be found in Robinson and Hidalgo
(\citeyear{robHid1997}).
Under Assumption~\ref{as3.1}, Robinson and Hidalgo (\citeyear{robHid1997}) showed that
%
\begin{equation}\label{e3.2}
\sqrt{n}(\hat{\alpha}_{n}-\alpha)\stackrel{\mathcal{L}}{\to}
N(0, \Sigma_{\phi}^{-1}\Sigma_{\psi}\Sigma_{\phi}^{-1}) .
\end{equation}
The intercept term $\alpha_{0}$ is estimated by
\[
\hat{\alpha}_{0n}=\bar{y}-\hat{\alpha}_{n}'\tilde{x}=\alpha_{0}
+\bar{\varepsilon}- (\hat{\alpha}_{n}-\alpha)'\bar{x} ,
\]
where
$\bar{\varepsilon}=\sum_{t=1}^{n}\varepsilon_{t}/n$. When $H\in(1/2, 1)$
or $H\le1/2$ with $Ex_{t}=0$, we see that
$n\sigma_{n}^{-1}(\hat{\alpha}_{n}-\alpha)'\bar{x}=o_{p}(1)$ and hence,
in these cases, we have
%
\begin{equation}\label{e3.3}
n\sigma_{n}^{-1}(\hat{\alpha}_{0n}-\alpha_{0})\stackrel{\mathcal{L}}{\to} N(0, 1).
\end{equation}

The results of Robinson and Hidalgo (\citeyear{robHid1997}) hold not only for
long-memory $\{\varepsilon_{t}\}$ but also for short-memory
$\{\varepsilon_{t}\}$. The following result entails the residual
empirical process for both long- and short-memory cases.
\begin{theorem}\label{thm3.1}
If Assumptions \ref{as2.1} and \ref{as3.1} hold, then the
results of Theorem~\ref{thm2.1} hold with
$\hat{\beta}_{n}=(\hat{\alpha}_{0n}, \hat{\alpha}_{n}')'$,
$\delta_{n}= \operatorname{diag}(\sigma_{n}n^{-1}, n^{-1/2}I_{q})$ and
$X_{t}=(1, x_{t}')'$.
\end{theorem}
\begin{pf}
It is readily seen that Assumptions \ref{as2.2}(a)--(c) hold. Note that
\begin{eqnarray*}
&&\frac{1}{\sigma_{n}}\sup_{x}
\Biggl\|\sum_{t=1}^{n}A_{t}(x)\delta_{n}'X_{t}\Biggr\|
\le \sup_{x}\Biggl\|\frac{1}{n}\sum_{t=1}^{n}A_{t}(x)\Biggr\|
+ \frac{1}{\sqrt{n}\sigma_{n}} \sup_{x}
\Biggl\|\sum_{t=1}^{n}A_{t}(x)x_{t}\Biggr\|.
\end{eqnarray*}
To check Assumption~\ref{as2.2}(d), we only need to show that
%
\begin{eqnarray}\label{e3.4}
\sup_{x}\frac{1}{\sqrt{n}\sigma_{n}} \sup_{x}
\Biggl\|\sum_{t=1}^{n}A_{t}(x)x_{t} \Biggr\|=o_{p}(1).
\end{eqnarray}
Similarly, it can be proved that $\sup_{x}|\sum
_{t=1}^{n}A_{t}(x)|/n=o_{p}(1)$. Since\break $\sup_{x}|G_{0}''(x)| <\infty$
implies $\lim_{|x|\to\infty}G_{0}'(x)=0$ [see Lee and Wei (\citeyear{leeWei1999})],
we see that $E\sup_{|x|>M}\{G'_{0}(x-\xi_{t-1})\|x_{t}\|\}\to0$ as $M\to\infty$.
Since $\sqrt{n}/\sigma_{n}=O(1)$, for
any given $\epsilon>0$, there exists a constant $M>0$ such that
%
\begin{eqnarray}\label{e3.5}
&& P\Biggl(\sup_{|x|> M}\Biggl\|\frac{1}{\sqrt{n}\sigma_{n}}
\sum_{t=1}^{n}A_{t} (x)x_{t}\biggr\|>\eta\Biggr)
\nonumber\\[-8pt]
\\[-8pt]
&&\qquad \le\frac{2\sqrt{n}}{\sigma_{n}\eta}
E\sup_{|x|>M}\{G'_{0}(x-\xi_{t-1})\|x_{t}\|\} < \epsilon ,
\nonumber
\end{eqnarray}
uniformly in $n$. Partition $[-M, M]$ into $m=[4M\delta^{-1}]$
subintervals such that $-M=c_{0}\le c_{1}\le\cdots\le c_{m}=M$ with
$c_{r+1}-c_{r}<\delta$ for any given constant $\delta>0$. Let
$U_{nr}=(\sqrt{n}\sigma_{n})^{-1}\sum_{t=1}^{n} A_{t}(c_{r})x_{t}$.
When $H\in(1/2, 1)$, $\|U_{nr}\|\le2n^{-1/2-H}\times\break \sum_{t=1}^{n}
\|x_{t}\|=o_{p}(1)$. When $H\le1/2$, since $A_{t}(c_{r})$ and
$x_{t}$ are independent for each~$c_{r}$, we can show that
$U_{nr}=o_{p}(1)$. Thus, we have
%
\begin{eqnarray}\label{e3.6}
&& \sup_{|x|\le M} \Biggl\|\frac{1}{\sqrt{n}\sigma_{n}}
\sum_{t=1}^{n}A_{t}(x)x_{t}  \Biggr\|
\nonumber\\
&&\qquad \le \max_{r} \sup_{x\in[c_{r},c_{r+1}]}
\Biggl\|\frac{1}{\sqrt{n}\sigma_{n}}\sum_{t=1}^{n}
[A_{t}(x)-A_{t}(c_{r})]x_{t}\Biggr\|
 + \max_{r}\|U_{nr}\|
\nonumber\\[-8pt]
\\[-8pt]
&&\qquad \le 2\delta\sup_{x}|G''_{0}(x)|O_{p}(1)+o_{p}(1)
\nonumber\\
&&\qquad =O_{p}(\delta)+o_{p}(1) .
\nonumber
\end{eqnarray}
Using (\ref{e3.5})--(\ref{e3.6}), (\ref{e3.4}) is established.
\end{pf}

We see that $R_{n}=O_{p}(1)$ and $K_{n}(x)=O_{p}(1)$. When
$Ex_{t}=0$, we have
$R_{n}(x)=n\sigma_{n}^{-1}(\hat{\alpha}_{0n}-\alpha_{0}) \neq o_p(1)$
by virtue of (\ref{e3.3}). In this case,
the estimated mean affects the limit distribution of
$K_{n}(x)$ by Theorem~\ref{thm3.1}. By (\ref{e1.3}) and (\ref{e3.3}),
we have the following result.
\begin{corollary}\label{cor3.1}
If Assumptions \ref{as2.1} and \ref{as3.1}
hold and $H\in(1/2,1)$, then
\[
\biggl[
\sup_{x}F'(x)\biggr]^{-1}\sup_{ x }|\hat{K}_{n}(
x)|\stackrel{\mathcal{L}}{\rightarrow} |N(0,4)|.
\]
\end{corollary}
\begin{remark}\label{rem3.1}
This corollary gives a statistic for
testing the distribution of the long-memory noises in model (\ref{e3.1})
when $\alpha_0$ is unknown. The asymptotic variance of this test
statistic is four times bigger than that in (\ref{e1.5}), which reflects
the effects of the slower convergence rate of the estimated
parameter $\hat{\alpha}_{0n}$. When $\alpha_0$ is known, the test
statistic (\ref{e1.5}) is still valid, however. As pointed out by the
reviewer, when $F=F(x, \theta)$ involves an unknown parameter
$\theta$, one should consider $\hat{K}_{n}$ with $F(x)$ being
replaced by $F(x, \hat{\theta}_{n})$.
Under such circumstances, the limit
distribution of the statistic is usually different from that of
Corollary~\ref{cor3.1}. This fact serves as a reminiscence of the classical
Kolmogorov--Smirnov statistics problem when the underlying parameters
are estimated; see Durbin (\citeyear{dur1976}). When $H \le 1/2$, it can be shown
that the limit distribution of the statistic exists by means of the
result of Wu (\citeyear{wu2003}). The closed form of such a limit distribution is
rather complicated and does not possess a simple expression,
however, and is not presented here.
\end{remark}

\section{Residual empirical process of unstable AR($p$) models}\label{s4}

This section considers the unstable AR$(p)$ model with starting
value $\{y_{0}, y_{-1},\ldots, y_{-p+1}\}$ independent of
$\{\varepsilon_{s}\dvtx s<0\}$ such that
%
\begin{equation}\label{e4.1}
y_{t}=\beta' X_{t}+ \varepsilon_{t},
\end{equation}
where $X_{t}=(y_{t-1},\ldots, y_{t-p})'$, $\beta=
(\phi_{1},\ldots,\phi_{p})'$, and the characteristic polynomial
$\phi(z)=1-\phi_{1}z-\cdots-\phi_{p}z^{p} $ has the decomposition,
%
\begin{equation}\label{e4.2}
\phi(z)=(1-z)^{a}(1+z)^{b}\prod_{k=1}^{l}[(1-ze^{i\theta_{k}})(1
+ze^{i\theta_{k}})]^{d_{k}},
\end{equation}
$a,  b, l,  d_{k}, k=1,\ldots,l$, are nonnegative
integers, $p=a+b+2(d_{1}+\cdots+d_{l})$,
and $\{\varepsilon_{t}\}$ is defined in model (\ref{e1.1}). Here, $a$ denotes
the multiplicity of the root $z=1$ for $\phi(z)=0$. Same interpretations
are given to $b$ and $l$. We estimate $\beta$ by the LSE:
\[
\hat{\beta}_{n}=\Biggl(\sum_{t=1}^{n}X_{t}X_{t}'\Biggr)^{-1}
\sum_{t=1}^{n}X_{t}y_{t}.
\]

For the special case with $\phi(z)=1-z$, Wu (\citeyear{wu2006}) obtained the
limiting distribution of $\hat{\beta}_{n}$
under Assumption~\ref{as2.1}(a); see also Sowell (\citeyear{sow1990}) and
Wang, Lin and Gulati (\citeyear{wanLinGul2003}). For the general case,
the limit distribution of $\hat{\beta}_{n}$ was obtained by
Chan and Terrin (\citeyear{chaTer1995}) and
Truong-Van and Larramendy (\citeyear{troLar1996}) under
the following Assumption~\ref{as4.1}(a) and (b),
respectively. It can be seen that Assumption~\ref{as2.1}(a) is much weaker than
Assumption~\ref{as4.1}.
\begin{assumption}\label{as4.1}
(a) $L_{0}(j)\sim c$, $c$ is a constant, $H\in(1/2, 1)$ and
$e_{t}\sim N(0, \sigma^{2}_{e})$, or (b)
$\sum_{j=0}^{\infty}j|a_{j}|<\infty$ and $\sigma>0$.
\end{assumption}

Let $\delta_{n}=G'J_{n}^{-1}$, where $G$ is the constant matrix
given in Chan and Wei (\citeyear{chaWei1988}) and
$J_{n}=\operatorname{diag}(N_{1}, N_{2},\ldots, N_{l+2})$ with
$N_{1}=\operatorname{diag}(n, n^{2}, \ldots, n^{a})$,
$N_{2}=\operatorname{diag}(n, n^{2}, \ldots, n^{b})$ and
$N_{k+2}=\operatorname{diag}(nI_{2}, \ldots, n^{d_{k}}I_{2})$,
$k=1, \ldots,l$. Define\break $\xi_{H}(\tau)= [f_{0}(\tau),\ldots, f_{a-1}(\tau)]'$,
$f_{0}(\tau)=B_{H}(\tau)$ and
$f_{j}(\tau)=\int_{0}^{\tau}f_{j-1}(s) \,ds$, $j=1, \ldots, a$,
where $B_{H}(\tau)$ is a fractional Brownian motion with covariances
\[
E[B_{H}(\tau)B_{H}(s)]=\tfrac{1}{2}\{s^{2H}+\tau^{2H}-
|s-\tau|^{2H}\}  \qquad \mbox{for }  0\le s,\tau\le1.
\]
We now state
the results for model (\ref{e4.1}).
\begin{theorem}\label{thm4.1}
For model (\ref{e4.1}), if Assumption \ref{as2.1} holds
with $\phi(z)=1-z$, or if Assumption~\ref{as4.1}\textup{(a)} holds, or
if Assumptions \ref{as2.1}\textup{(b)} and \ref{as4.1}\textup{(b)} hold,
then the result of Theorem~\ref{thm2.1}
holds with $R_{n}=o_{p}(1)$ for
$a=0$ and
\[
R_{n} \stackrel{\mathcal{L}}{\longrightarrow}
\cases{%
\displaystyle
(\Gamma+\zeta_{1/2})'\Omega^{-1}_{1/2} \int_{0}^{1}
\xi_{1/2}(\tau)  \, d\tau, & \quad if $H\le1/2$,\cr
\displaystyle
\zeta'_{H}\Omega^{-1}_{H} \int_{0}^{1} \xi_{H}(\tau)
\,d\tau, & \quad if $H\in(1/2, 1)$,}
\]
for $a\ge1$, where $\Gamma=(\gamma, 0,\ldots, 0)'_{a\times1}$,
$\gamma=1/2(1-E\varepsilon_{t}^{2}/\sigma^{2})$,
$\zeta_{H}=\break \int_{0}^{1}\xi_{H}(\tau) \,dB_{H}(\tau)$,
$\Omega_{H}=(\omega_{ij})_{a\times a}$ and
$\omega_{ij}=\int_{0}^{1}f_{i}(\tau)f_{j}(\tau)\, d\tau$.
\end{theorem}

Let
$D[0,1]$ be the Skorokhod space and $D^{p}=D\times D \times
\cdots\times D$ denote the $p$-Cartesian product space of
$D=D[0,1]$. To prove Theorem~\ref{thm4.1}, we need the following lemma.
Using the results in Chan and Wei (\citeyear{chaWei1988}),
Truong and Larramendy (\citeyear{troLar1996}) and
Wu (\citeyear{wu2006}), its proof is similar to that of Lemma 2.1 in
Ling (\citeyear{lin1998}) and the details are omitted.
\begin{lemma}\label{lem4.1}
Let $\tilde{\xi}=\xi_{H}$ if $H\in(1/2,1)$ and
$\tilde{\xi}=\xi_{1/2}$ if $H\le1/2$. If the assumptions
of Theorem~\ref{thm4.1} hold, then:
\begin{eqnarray*}
&& \mbox{\textup{(a)}} \quad \frac{1}{\sigma_{n}}\sum_{t=1}^{[n\tau]}\delta_{n}'X_{t}
\stackrel{\mathcal{L}}{\longrightarrow}
\biggl(\int_{0}^{\tau} \tilde{\xi}'(s)  \,ds, 0\biggr)'
\qquad  \mbox{in }   D^{p},   \mbox{ if }  a\ge1,
\\
&& \mbox{\textup{(b)}} \quad
\frac{1}{\sigma_{n}}\sum_{t=1}^{[n\tau]}\delta_{n}'X_{t}=o_{p}(1)
\qquad \mbox{uniformly for all }   \tau\in[0,1]  \mbox{ if }  a=0,
\\
&& \mbox{\textup{(c)}}\quad  \frac{1}{\sigma_{n}}\sum_{t=1}^{n}E\|\delta_{n}'X_{t}\|=O(1),
\\
&& \mbox{\textup{(d)}} \quad \frac{n}{\sigma_{n}^{2}}
\sum_{t=1}^{n}E\|\delta_{n}'X_{t}\|^{2}=O(1).
\end{eqnarray*}
\end{lemma}
\begin{pf}
For simplicity, we only prove Theorem~\ref{thm4.1}
for $\phi(z)=(1-z)$, that is, model (\ref{e4.1}) only has one unit root. The
general case can similarly be proved by Lemma \ref{lem4.1}.
When $\phi(z)=(1-z)$, $\delta_{n}=n^{-1}$ and
$X_{t}=y_{t-1}=\sum_{i=1}^{t-1}\varepsilon_{i}$. By Theorem~6.1 of
Chan and Terrin (\citeyear{chaTer1995}) and Theorem~3.1 of Truong-Van and Larramendy
(\citeyear{troLar1996}) or Theorems 3 and 4 of Wu (\citeyear{wu2006}),
Assumption~\ref{as2.2}(a) holds. By
Lemma \ref{lem4.1}(c) and (d), we see that Assumption~\ref{as2.2}(b)~and~(c) holds.

We now consider Assumption~\ref{as2.2}(d). First, note that
$E\sup_{|x|>M}A_{t}^{2}(x)\to0$ as $M\to\infty$ and
$\max_{1\le t\le n}\sigma_{n}^{-2}EX_{t}^{2}=O(1)$. Thus, for any
given $\epsilon>0$ and $\eta>0$, there exists a constant $M>0$ such
that
%
\begin{eqnarray}\label{e4.3}
&& P\Biggl(\sup_{|x|>M} \Biggl|\frac{1}{n\sigma_{n}}
\sum_{t=1}^{n}A_{t}(x)X_{t} \Biggr|>\eta\Biggr)
\nonumber\\[-8pt]
\\[-8pt]
&&\qquad \le \frac{\sqrt{E \sup_{|x|>M}|A_{t}(x)|^{2}}}{\eta n\sigma_{n}}
\sum_{t=1}^{n}\sqrt{E|X_{t}|^{2}} < \epsilon,
\nonumber
\end{eqnarray}
uniformly in $n$. Partition $[-M, M]$ into $m=[4M\delta^{-1}]$
subintervals such that $-M=x_{0}\le x_{1}\le\cdots\le x_{m}=M$ with
$x_{r+1}-x_{r}<\delta$ for any given $\delta>0$. Thus,
%
\begin{eqnarray}\label{e4.4}
&& \sup_{|x|\le M} \Biggl|\frac{1}{n\sigma_{n}}
\sum_{t=1}^{n}A_{t}(x)X_{t} \Biggr|
\nonumber\\
&&\qquad \le \max_{ r} \sup_{x_{r-1}\le x\le x_{r}}
\Biggl |\frac{1}{n\sigma_{n}}\sum_{t=1}^{n}A_{t}(x)X_{t} \Biggr|
\nonumber\\[-8pt]
\\[-8pt]
&&\qquad \le \max_{ r} \sup_{x_{r-1}\le x\le x_{r}}
\Biggl |\frac{1}{n\sigma_{n}}\sum_{t=1}^{n}
[A_{t}(x)-A_{t}(x_{r})]X_{t} \Biggr|
\nonumber\\
&&\qquad \quad {} + \max_{r}
\Biggl |\frac{1}{n\sigma_{n}}\sum_{t=1}^{n}A_{t}(x_{r})X_{t} \Biggr|
= J_{1n}+J_{2n} , \qquad \mbox{say}.
\nonumber
\end{eqnarray}
Since $\sup_{x}|A_{t}'(x)|<\infty$, by Lemma \ref{lem4.1}(c) and
the Taylor expansion, we have
%
\begin{equation}\label{e4.5}
J_{1n}\le O(\delta) \Biggl[\frac{1}{n\sigma_{n}}\sum_{t=1}^{n}|X_{t}|
 \Biggr]=O_{p}(\delta).
\end{equation}

For $J_{2n}$, we need the following decomposition:
\begin{eqnarray*}
\frac{1}{n\sigma_{n}}\sum_{t=1}^{n}A_{t}(x)X_{t}
&=& \frac{1}{n\sigma_{n}}\sum_{i=1}^{n}
\Biggl[\sum_{t=i+1}^{n}A_{t}(x)\Biggr]\varepsilon_{i}
\\
&=& \frac{1}{n\sigma_{n}}\Biggl(\sum_{i=1}^{n}\varepsilon_{i}\Biggr)
\Biggl[\sum_{t=1}^{n}A_{t}(x)\Biggr]
- \frac{1}{n\sigma_{n}}\sum_{i=1}^{n}
\Biggl[\sum_{t=1}^{i}A_{t}(x)\Biggr]\varepsilon_{i}
\\
&=&U_{1n}(x)-U_{2n}(x),  \qquad \mbox{say}.
\end{eqnarray*}
By the ergodic theorem, $\sum_{t=1}^{n}A_{t}(x)/n=o_{p}(1)$ for each
$x$. Furthermore,
since $\sum_{i=1}^{n}\varepsilon_{i}/\sigma_{n}=O_{p}(1)$, we have
$\max_{r}|U_{1n}(x_{r})|=o_{p}(1)$ for a given $\delta>0$.

We next consider $U_{2n}(x)$. When $H\le1/2$, by Theorem 2 of Wu
(\citeyear{wu2006}), we know that $\sum_{t=1}^{[n\tau]}A_{t}(x)/\sigma_{n}
\stackrel{\mathcal{L}}{\to} S(\tau)$ in~$D$ for each $x$ and
$\sum_{t=1}^{[n\tau]}\varepsilon_{t}/\sqrt{n}\stackrel{\mathcal{L}}{\to}
\xi(\tau)$ in~$D$, where $S(\tau)$ and $\xi(\tau)$ are standard Brownian
motions. By Theorem~\ref{thm3.1} of Ling and Li (\citeyear{linLi1998}),
$U_{2n}(x)=o_{p}(1)$
for each $x$ and, hence, $\max_{r}|U_{2n}(x_{r})|=o_{p}(1)$ for any
given $\delta>0$. Thus, Assumption~\ref{as2.2}(d) holds when $H\le1/2$.

When $H\in(1/2, 1)$, we decompose $U_{2n}(x)$ as follows:
%
\begin{equation}\label{e4.6}
\qquad
\frac{1}{n\sigma_{n}}\sum_{i=1}^{n}
\Biggl[\sum_{t=1}^{i}R_{t}(x)\Biggr]\varepsilon_{i}
+ \frac{G''_{0}(x)}{n\sigma_{n}}\sum_{i=1}^{n}
\Biggl(\sum_{t=1}^{i}\xi_{t-1}\Biggr)\varepsilon_{i}
= U_{3n}(x)+U_{4n}(x),
\end{equation}
say, where $R_{t}(x)= A_{t}(x)-G''_{0}(x)\xi_{t-1}$. For each $x$
and any $\zeta>0$, by Corollary 1 of Wu (\citeyear{wu2006})
[see also Theorem~\ref{thm3.1} in Ho and Hsing (\citeyear{hoHsi1997})], we have
%
\begin{equation}\label{e4.7}
E\Biggl[\sum_{t=1}^{i}R_{t}(x)\Biggr]^{2}
= O\bigl(i^{\max\{1, 4(H-1/2)+2\zeta\}}\bigr).
\end{equation}
By (\ref{e4.7}), for any $\eta>0$ and $\delta>0$, we have
%
\begin{eqnarray}\label{e4.8}
P\biggl(\max_{r}|U_{3n}(x_{r})|>\eta\biggr)
&\le& \frac{1}{\eta}\sum_{r=1}^{m}E|U_{3n}(x_{r})|
\nonumber\\
&\le& \frac{1}{\eta n \sigma_{n}}\sum_{r=1}^{m}\sum_{i=1}^{n}
\Biggl\{E\Biggl[\sum_{t=1}^{i} R_{t}(x)\Biggr]^{2}E\varepsilon_{i}^{2} \Biggr\}^{1/2}
\\
& =& O(n^{-\gamma}L_{0}^{-1}(n))\to0,
\nonumber
\end{eqnarray}
when $n \to\infty$, where $\gamma=\min\{H-1/2,  1-H-\zeta\}>0$.
Note that
\begin{eqnarray*}
U_{4n}(x)&=&- \frac{G''_{0}(x)}{n\sigma_{n}}\sum_{i=1}^{n}
\Biggl(\sum_{t=1}^{i}\varepsilon_{t} \Biggr)\varepsilon_{i}
+ \frac{G''_{0}(x)}{n\sigma_{n}}\sum_{i=1}^{n}
\Biggl(\sum_{t=1}^{i}e_{t} \Biggr)\varepsilon_{i} .
\end{eqnarray*}
By Theorems 3.2 and 3.3 of Chan and Terrin (\citeyear{chaTer1995})
or Theorem 3 of Wu (\citeyear{wu2006}),
\[
\sum_{i=1}^{n} \Biggl(\frac{1}{\sigma_{n}}
\sum_{t=1}^{i}\varepsilon_{t} \Biggr)\frac{\varepsilon_{i}}{\sigma
_{n}}\stackrel{\mathcal{L}}{\longrightarrow} \int_{0}^{1}B_{H}(s) \,dB_{H}(s).
\]
Thus, the first term in $U_{4n}(x)$ is
$o_{p}(1)$ uniformly in $x\in{R}$. Note that\break
$\sum_{t=1}^{n}|\varepsilon_{t}|/n=O_{p}(1)$ by the ergodic theorem
and $\max_{1\le i\le
n}|\sum_{t=1}^{i}e_{t}|/\sqrt{n}\stackrel{\mathcal{L}}{\to}
\max_{0\le\tau\le1}|B_{1/2}(\tau)|$. Since
$\sqrt{n}/\sigma_{n}=O(n^{-H+1/2}/L_{0}(n))=o(1)$, the second term
in $U_{4n}(x)$ is $o_{p}(1)$ uniformly in $x\in{R}$. Thus, we have
$\max_{ x}|U_{4n}(x)|=o_{p}(1)$.
Furthermore, by (\ref{e4.6}) and (\ref{e4.8}), $\max_{ r}|U_{2n}(x_{r})|=o_{p}(1)$ for
any given $\delta$ when $H\in(1/2, 1)$. Thus, Assumption~\ref{as2.2}(d)
holds when $H\in(1/2,1)$.
\end{pf}
\begin{remark}\label{rem4.1}
From this theorem, we see that the
empirical process of $\{{\varepsilon}_{t}\}$ is not affected if
$\{{\varepsilon}_{t}\}$ is replaced by $\{\hat{\varepsilon}_{t}\}$
when $\phi(z)$ does not have a root equaling one. It has a
profound effect when $\phi(z)$ has a unit root, however. In
particular, using Theorem 3 of Wu (\citeyear{wu2006}), we have the following
corollary.
\end{remark}
\begin{corollary}\label{cor4.1}
If $\phi(z)=(1-z)$ and Assumption~\ref{as2.1} holds with
$H\in(1/2,1)$, then it follows that
\begin{eqnarray*}
&& \biggl[
\sup_{x}F'(x)\biggr]^{-1}
\sup_{ x}|\hat{K}_{n}(x)|
\\
&&\qquad\stackrel{\mathcal{L}}{\longrightarrow}
\biggl |B_{H}(1)+\biggl[\int_{0}^{1}B_{H}(\tau)\,dB_{H}(\tau)\biggr]
\biggl[\int_{0}^{1}B_{H}(\tau) \,d\tau\biggr]
\biggl[\int_{0}^{1}B_{H}^{2}(\tau)\, d\tau\biggr]^{-1} \biggr|.
\end{eqnarray*}
\end{corollary}
\begin{remark}\label{rem4.2}
Corollary~\ref{cor4.1} gives the limit
distribution of the Kolmogorov--Smirnov statistic. It can be
used to test for the distribution of the long-memory noises in
model (\ref{e4.1}). For instance, using $\hat{\varepsilon_{t}}$ as a proxy
for $\varepsilon_{t}$, $H$ may be estimated by Robinson's (\citeyear{rob1995a})
semiparametric method. Although the asymptotic validity of such a
procedure still needs to be examined, for a given $H\in(1/2, 1)$,
the percentiles of the limit distribution can be tabulated by means
of simulations. Corollary~\ref{cor4.1} thus provides a means to apply the
Kolmogorov--Smirnov statistics to model (\ref{e4.1}).
\end{remark}

\begin{appendix}
\setcounter{lemma}{0}
\setcounter{equation}{0}
\section*{Appendix: Technical lemmas}

Let $x_{r}=r\epsilon\sigma_{n}^{-1}$ for any $r\in Z$ and
some $\epsilon>0$ and decompose the real line $R$ as
$R=\bigcup_{r\in Z}[x_{r}, x_{r+1}]$.
Let $g_{t}(u, \lambda)$ be defined in (\ref{e2.3}) and
\[
a_{nt}(x)=I\bigl(\varepsilon_{t}\le x+g_{t}(u,\lambda)\bigr)
- F_{t-1}(x) -I(\varepsilon_{t}\le x)+G_{0}(x-\xi_{t-1}),
\]
where $ F_{t-1}(x)=
E[I(e_{t}\le x-\xi_{t-1}+g_{t}(u, \lambda))|\mathcal{F}_{t-1}]=G_{0}[
x-\xi_{t-1}+g_{t}(u, \lambda)]$, $u\in[-\Delta, \Delta]^{p}$ with
$\Delta>0$ and $\lambda\in[-1, 1]$. We have the following lemma.
\begin{lemma}\label{lemA.1}
Let $\tilde{Z}_{1n}(x, u,\lambda)=\sum_{t=1}^{n}a_{nt}(x)/\sigma_{n}$.
For every $u$ and $\lambda$, if Assumption~\ref{as2.1} and
Assumptions \ref{as2.2}\textup{(b)} and \textup{(c)} hold, then:
\begin{eqnarray*}
&&\mbox{\textup{(a)}}\quad
\max_{r}\max_{x\in[x_{r}, x_{r+1}]}\frac{1}{\sigma_{n}}\sum_{t=1}^{n}
\bigl|F\bigl(x_{r+1}+g_{t}(u, \lambda)\bigr)
- F\bigl(x+g_{t}(u,\lambda)\bigr)\bigr|=O_{p}(\epsilon),
\\
&& \mbox{\textup{(b)}}\quad
\sup_{r}|\tilde{Z}_{1n}(x_{r}, u,\lambda)|=o_{p}(1)  \qquad
\mbox{for any given $\epsilon>0$}.
\end{eqnarray*}
\end{lemma}
\begin{pf}
By Assumption~\ref{as2.1}(b), $F'(x)$ exists and is
bounded; see Ho and Hsing (\citeyear{hoHsi1996}).
Since $n/\sigma^{2}_{n}=O(1)$, by the
Taylor expansion, part (a) holds.

For part (b), since $\sum_{t=1}^{n}a_{nt}(x)$ is a martingale array with
respect to $\mathcal{F}_{n}= \sigma\{(e_{t}, X_{t})$, $ t\le n\}$, by
the Rosenthal inequality [see page 23 of Hall and Heyde (\citeyear{halHey1980})],
%
\begin{eqnarray}\label{eA.1}
E\Biggl[\sum_{t=1}^{n}a_{nt}(x)\Biggr]^{4}
&\le& c E\Biggl\{\sum_{t=1}^{n}E[a_{nt}^{2}(x)|\mathcal{F}_{t-1}]\Biggr\}^{2}
+ c\sum_{t=1}^{n}E[a_{nt}^{4}(x)]
\nonumber\\[-8pt]
\\[-8pt]
&\le& cn \sum_{t=1}^{n}E\{E[a_{nt}^{2}(x)|\mathcal{F}_{t-1}]\}^{2}
+ 2c\sum_{t=1}^{n}E[a_{nt}^{2}(x)]
\nonumber
\end{eqnarray}
for some constant $c$, where we use $a_{nt}^{4}(x)\le2a_{nt}^{2}(x)$.
Denote $g_{t}(u,\lambda)$ by $g_{t}$ and let $H_{t}^{\pm}(x)=G_{0}(x-\xi
_{t-1}\pm|g_{t}|)$. Since
$E[I( e_{t}\le x-\xi_{t-1})|\mathcal{F}_{t-1}]=G_{0}( x-\xi_{t-1})$ and
$G_{0}(x)$ is nondecreasing, we have
\begin{eqnarray*}
E[a_{nt}^{2}(x)|\mathcal{F}_{t-1} ] \le|F_{t-1}(x)
- G_{0}(x-\xi_{t-1})|\le H_{t}^{+}(x)-H_{t}^{-}(x).
\end{eqnarray*}
Again, since $G_{0}(x)$ is nondecreasing, for any positive integer
$M$, we have
%
\begin{eqnarray}\label{eA.2}
&& \sum_{r=-M}^{M}E[H_{t}^{+}(x_{r})-H_{t}^{-}(x_{r})]
\nonumber\\
&&\qquad \le  \frac{{\sigma_{n}}}{\epsilon}\sum_{r=-M}^{M}
E\biggl[\int_{x_{r}}^{x_{r+1}}  H_{t}^{+}(x) \,dx
- \int_{x_{r-1}}^{x_{r}}  H_{t}^{-}(x) \,dx\biggr]
\nonumber\\
&&\qquad =  \frac{{\sigma_{n}}}{\epsilon}E
\biggl\{\int_{x_{M}}^{x_{M+1}}H_{t}^{+}(x) \,dx
+ \int_{x_{-M-1}}^{x_{-M}}H_{t}^{-}(x) \,dx
\nonumber\\
&&\hspace*{98pt}{} +\int_{x_{-M}}^{x_{M}}
[H_{t}^{+}(x)-H_{t}^{-}(x)] \,dx \biggr\}
\\
&&\qquad \le  2+\frac{{\sigma_{n}}}{\epsilon}E  \biggl\{
\int_{x_{-M}}^{x_{M}} \int_{-|g_{t}|}^{|g_{t}|}
G_{0}'(x-\xi_{t-1}+y) \,dy \,dx \biggr\}
\nonumber\\
&&\qquad \le 2+\frac{{\sigma_{n}}}{\epsilon}E  \biggl\{
\int_{-|g_{t}|}^{|g_{t}|} \int_{-\infty}^{\infty}
G_{0}'(x-\xi_{t-1}+y)\, dx \,dy \biggr\}
\nonumber\\
&&\qquad = 2+\frac{2\sigma_{n}}{\epsilon}E |g_{t}|.
\nonumber
\end{eqnarray}
Similarly, we have
%
\begin{eqnarray}\label{eA.3}
\qquad
 \sum_{r=-M}^{M}E[H_{t}^{+}(x_{r})-H_{t}^{-}(x_{r})]^{2}
&\le&
c\sum_{r=-M}^{M}E\{|g_{t}|[H_{t}^{+}(x_{r})-H_{t}^{-}(x_{r})]\}
\nonumber\\[-8pt]
\\[-8pt]
&=& 2c E|g_{t}|+\frac{2c\sigma_{n}}{\epsilon}Eg_{t}^{2},
\nonumber
\end{eqnarray}
where $c=2\sup_{x}G_{0}'(x)$. Using (\ref{eA.2})--(\ref{eA.3})
and Assumptions \ref{as2.2}(b)--(c),
%
\begin{eqnarray}\label{eA.4}
\qquad
\frac{ 1}{\sigma_{n}^{4}}\sum_{r}\sum_{t=1}^{n}
E[a_{nt}^{2}(x_{r})]
&\le& \frac{1}{\sigma_{n}^{4}}\lim_{M\to\infty}\sum_{r=-M}^{M}
\sum_{t=1}^{n}E[H_{t}^{+}(x_{r})-H_{t}^{-}(x_{r})]
\nonumber\\[-8pt]
\\[-8pt]
& \le& \frac{ 2n}{\sigma_{n}^{4}}
+ \frac{2}{\epsilon\sigma_{n}^{3}}\sum_{t=1}^{n} E |g_{t}|
=o(1),
\nonumber
\end{eqnarray}
%
\begin{equation}\label{eA.5}
\qquad \quad
\frac{ n}{\sigma_{n}^{4}}\sum_{r}\sum_{t=1}^{n} E\{E
[a_{nt}^{2}(x_{r})|\mathcal{F}_{t-1}]\}^{2}
\le\frac{2n}{\sigma_{n}^{4}}\sum_{t=1}^{n} E |g_{t}|
+ \frac{2}{\epsilon\sigma_{n}^{3}}\sum_{t=1}^{n} E g_{t}^{2}=o(1),
\end{equation}
as $n/\sigma_{n}^{2}=O(1)$.
By the Markov inequality, (\ref{eA.1}), (\ref{eA.4}) and (\ref{eA.5}),
\begin{eqnarray*}
P\biggl(\sup_{r}|\tilde{Z}_{1n}(x_{r}, u,\lambda)|\ge\eta\biggr)
&\le& \sum_{r}P\bigl( |\tilde{Z}_{1n}(x_{r}, u,\lambda)|\ge\eta\bigr)
\\
&\le& \frac{1}{\eta^{4}\sigma_{n}^{4}}\sum_{r}E
\Biggl[\sum_{t=1}^{n}a_{nt}(x_{r})\Biggr]^{4}
\\
&=& o(1),
\end{eqnarray*}
as $n\to\infty$, for any given $\epsilon>0$. Thus, part (b) is proved.
\end{pf}
\begin{lemma}\label{lemA.2}
Let $\tilde{Z}_{2n}(x, u, \lambda)= \sum_{t=1}^{n}[F_{t-1}(x)-G_{0}(
x-\xi_{t-1})-F(x+g_{t}(u,\lambda))+F(x)]/\sigma_{n}$. If Assumptions
\ref{as2.1} and \ref{as2.2}\textup{(b)--(d)} hold, then\break
$\tilde{Z}_{2n}(x, u,\lambda)=\lambda J_{1n}(x)+J_{2n}(x, u,
\lambda)$ such that $\sup_{x}|J_{1n}(x)|=O_{p}(1)$ and\break
$\sup_{x}\sup_{u}\sup_{\lambda}|J_{2n}(u, x, \lambda)|=o_{p}(1)$.
\end{lemma}
\begin{pf}
By Assumption~\ref{as2.1}(b) and Lemma 6.2 of Ho and Hsing
(\citeyear{hoHsi1996}), $F''(x)$ exists
and is bounded. By the Taylor expansion and Assumption~\ref{as2.2}(c),
\begin{eqnarray*}
\tilde{Z}_{2n}(x, u, \lambda)
&=& \frac{1}{\sigma_{n}}\sum_{t=1}^{n}
\biggl\{A_{t}(x)g_{t}(u, \lambda)
+\frac{1}{2}g_{t}^{2}(u, \lambda)[G^{\prime\prime}_{0}(\xi_{t-1}^{*})
- F^{\prime\prime}(\tilde{\xi}_{t-1}^{*})]\biggr\}
\\
&=& \frac{1}{ \sigma_{n}}\sum_{t=1}^{n}A_{t}(x)g_{t}(u,\lambda)+o_{p}(1)
\\
&=& \frac{\lambda}{\sigma_{n}}\sum_{t=1}^{n}A_{t}(x)\|\delta_{n}'X_{t}\|
+ \Biggl[\frac{u}{\sigma_{n}}\sum_{t=1}^{n}A_{t}(x)\delta_{n}'X_{t}+o_{p}(1) \Biggr]
\\
&=&\lambda J_{1n}(x)+J_{2n}(x, u, \lambda) ,\qquad  \mbox{say},
\end{eqnarray*}
where we use $F'(x)=EG'_{0}(x-\xi_{t-1})$,
$\xi_{t-1}^{*}=x-\xi_{t-1}+\theta g_{t}(u, \lambda)$ and $\tilde{\xi
}_{t-1}^{*}=x+\tilde{\theta} g_{t}(u, \lambda)$ with
$\theta, \tilde{\theta}\in(0,1)$
and $o_{p}(1)$ being held uniformly in $x,u,\lambda$.
Since $\sup_{x}|A_{t}(x)|\le2$, by Assumption~\ref{as2.2}(b),
$\sup_{x}|J_{1n}(x)|=O_{p}(1)$.
Since $u\in[-\Delta, \Delta]^{p}$, by Assumption~\ref{as2.2}(d),
$\sup_{x}\sup_{u}\sup_{\lambda}|J_{2n}(x, u, \lambda)|=o_{p}(1)$.
The desired conclusion follows.
\end{pf}
\begin{lemma}\label{lemA.3}
If Assumptions \ref{as2.1} and \ref{as2.2}\textup{(b)--(d)}
hold, then it follows that
\begin{eqnarray*}
&&\sup_{x}|\tilde{Z}_{n}(x, u, \lambda)|\le J_{3n}(u, \lambda)
+ |\lambda|J_{4n},
\end{eqnarray*}
where $\tilde{Z}_{n}(x, u, \lambda)$ is
defined in (\ref{e2.3}), $0<J_{3n}(u, \lambda)=o_{p}(1)$ for each $u$ and
$\lambda$, and $0<J_{4n}=O_{p}(1)$ is independent of $u$.
\end{lemma}
\begin{pf}
Since $I(\varepsilon_{t}\le x)$ and
$F(x)$ are nondecreasing, for any $x\in[x_{r}, x_{r+1}]$,
\begin{eqnarray*}
\tilde{Z}_{n}(x, u, \lambda)
& \le& \tilde{Z}_{n}(x_{r+1}, u, \lambda)
+ \frac{1}{\sigma_{n}}\sum_{t=1}^{n}[F(x_{r+1}+g_{t})-F(x+g_{t})]
\\
&&{} + \frac{1}{\sigma_{n}}\sum_{t=1}^{n}
[I(\varepsilon_{t}\le x_{r+1})-F(x_{r+1})
- I(\varepsilon_{t}\le x)+F(x)] ,
\end{eqnarray*}
where $g_{t}$ denotes $g_{t}(u, \lambda)$ and a reverse inequality
holds when $x_{r+1}$ is replaced by $x_{r}$. Since
$|\tilde{Z}_{n}(x_{r+1}, u, \lambda)|\le |\tilde{Z}_{1n}(x_{r+1}, u,
\lambda)|+|\tilde{Z}_{2n}(x_{r+1}, u, \lambda)|$, we have
\[
\sup_{x}|\tilde{Z}_{n}(x ,u, \lambda)|
\le \max_{r}| \tilde{Z}_{2n}(x_{r}, u, \lambda)|+R_{n}(u, \lambda),
\]
where
%
\begin{eqnarray}\label{eA.6}
\qquad
R_{n}(u, \lambda)&=&\max_{r}| \tilde{Z}_{1n}(x_{r}, u, \lambda)|
\nonumber\\
&&{} + \max_{r}\max_{x\in[x_{r}, x_{r+1}]}\frac{1}
{\sigma_{n}}\sum_{t=1}^{n}|F(x_{r+1}+g_{t})-F(x+g_{t})|
\nonumber\\[-8pt]
\\[-8pt]
&&{} + \sup_{ |x_{1}-x_{2}|\le\epsilon\sigma^{-1}_{n}}\frac{1}{\sigma_{n}}
\Biggl |\sum_{t=1}^{n}[I(\varepsilon_{t}\le x_{1})
\nonumber\\
&&\hspace*{99pt}{}
- F(x_{1})-I(\varepsilon_{t}\le x_{2})+F(x_{2})] \Biggr|.
\nonumber
\end{eqnarray}
For any $\varepsilon, \eta>0$, by Lemma \ref{lem4.1}(a), we can take
$\epsilon$ small enough such that the second term of (\ref{eA.6}) is less
than $\eta$ happens with probability being at least $1-\varepsilon/4$.
For this $\epsilon$, the
first term of~(\ref{eA.6}) is $o_{p}(1)$ by Lemmas \ref{lemA.1}(b), and the
last term of~(\ref{eA.6}) is $o_{p}(1)$ by the tightness of the empirical
process of
$\{\varepsilon_{t}\}$ of Ho and Hsing (\citeyear{hoHsi1996})
and Wu (\citeyear{wu2003}). Thus,
$R_{n}(u, \lambda)=o_{p}(1)$ for each $u$ and $\lambda$. By virtue of
Lemma \ref{lemA.2}, the conclusion holds.
\end{pf}

\end{appendix}

\section*{Acknowledgments}

The authors would like to thank two
referees, an Associate Editor and the Co-editor, Professor M. L. Eaton,
for their helpful and constructive comments, which substantially
improved the presentation of this paper.

\printaddresses

\end{document}